\documentclass[11pt]{amsart}

\title{Injective convolution operators on $\lp{\infty}(\Gm)$ are surjective}
\author {Yemon Choi}
\address{Department of Mathematics, University of Manitoba, Winnipeg, R3T 2N2, Canada}
\email{y.choi.97@cantab.net}
\usepackage{amsmath,amssymb,amsfonts}
\usepackage{amsthm}   %Inbuilt environments for defns etc
\usepackage{concrete}

\theoremstyle{plain}
\newtheorem{thm}{Theorem}
 \newtheorem{propn}[thm]{Proposition}
\newtheorem{lemma}[thm]{Lemma}
\theoremstyle{definition}

\newtheorem{eg}[thm]{Example}
\theoremstyle{remark}
\newtheorem*{rem}{Remark}

\usepackage{amsmath,amssymb,amsfonts}
\newcommand{\dt}[1]{{\sf #1}}
\newcommand{\id}[1][]{{\mathbf 1}_{#1}}

\newcommand{\Cplx}{{\mathbb C}}
\newcommand{\Real}{{\mathbb R}}
\newcommand{\Torus}{{\mathbb T}}
\newcommand{\Zahl}{{\mathbb Z}}
\newcommand{\lp}[1]{\ell^{#1}}

\newcommand{\Gm}{\Gamma}
\newcommand{\om}{{\omega}}
\newcommand{\abs}[1]{\vert #1 \vert}

\newcommand{\FC}[1]{{\widehat{#1}}}

\begin{document}
\begin{abstract}
Let $\Gm$ be a discrete group and let $f \in \lp{1}(\Gm)$. We observe that if the natural convolution operator $\rho_f:\lp{\infty}(\Gm)\to \lp{\infty}(\Gm)$ is injective, then $f$ is invertible in $\lp{1}(\Gm)$. Our proof simplifies and generalizes calculations in a preprint of Deninger and Schmidt, by appealing to the direct finiteness of the algebra $\lp{1}(\Gm)$.

We give simple examples to show that in general one cannot replace $\lp{\infty}$ with $\lp{p}$, $1\leq p< \infty$, nor with $L^{\infty}(G)$ for nondiscrete $G$. Finally, we consider the problem of extending the main result to the case of weighted convolution operators on $\Gm$, and give some partial results.
\end{abstract}
\maketitle

\begin{section}{Introduction}
We begin with some background to motivate the result announced in the abstract.

Given a discrete group $\Gm$ and an element $f$ in its integral group ring $\Zahl\Gm$, one can construct a dynamical system with interesting properties (a so-called \dt{algebraic $\Gm$-action} on the Pontryagin dual of the additive group $\Zahl\Gm/\Zahl\Gm f$, see \cite{Den_JAMS06} for further details). This construction is studied in \cite{Den_JAMS06,DenSch}, where entropy formulas are given under the further assumption that $\Gm$ is residually finite and amenable.

In analyzing this action (e.g.~finding when it is expansive) one is led to consider the `convolution operator' $\rho_f: \lp{\infty}(\Gm)\to\lp{\infty}(\Gm)$ defined by
\begin{equation}\label{eq:rho-def}
\rho_f(g)(x)= \sum_{y\in\Gm} g(xy)f(y) \qquad\text{for all $g\in\lp{\infty}(\Gm)$ and all $x\in \Gm$.}
\end{equation}
Note that this definition makes sense for any $f\in\lp{1}(\Gm)$.
An early preprint version of \cite{DenSch} contains the following technical result.

\medskip
\noindent{\bf Theorem A }\cite[Theorem~1.4]{DenSch_pre} Let $\Gm$ be a countable and residually finite group and let $f\in\lp{1}(\Gm)$. If $\rho_f$ is injective then it is surjective.
\medskip

The original proof used the hypotheses on $\Gm$ to make various approximation arguments, exploiting the fact that injective operators on finite-dimensional vector spaces are surjective. G.~Elek (unpublished note) extended these techniques to prove Theorem A for all discrete groups $\Gm$ that are \dt{sofic}: we shall not give the definition here, but merely observe that
\begin{itemize}
\item there is no known example of a discrete group that is not sofic;
\item there is as yet no proof that all groups are sofic.
\end{itemize}

In this short note we show (Theorem~\ref{t:mainthm} below) that Theorem A holds for \emph{any} discrete group $\Gm$, by combining a simple duality argument with the old result that all left invertible elements in the group algebra $\lp{1}(\Gm)$ are automatically invertible.

One might naively hope to generalize Theorem~\ref{t:mainthm} in two directions: first, we could restrict $\rho_f$ to $c_0(\Gm)$ or $\lp{p}(\Gm)$ for $1\leq p <\infty$; and secondly we could study the analogous convolution operators on $L^\infty(G)$ for a (non-discrete) locally compact group $G$. We give simple examples to show that in both cases injectivity of $\rho_f$ need not imply surjectivity, even if we work with abelian groups.

In the last section we state some partial results for the more general case of \emph{weighted} convolution operators on $\lp{\infty}(\Gm)$; these cover the case of an arbitrary weight on an amenable group, which may be of interest given the original motivating examples in \cite{DenSch}, and also the case of symmetric weights on arbitrary groups.
We are unable to resolve the case of arbitrary weights on arbitrary groups, since it appears to be unknown if there is a discrete group $\Gm$ and a weight $\om$ on it such that the algebra $\lp{1}(\Gm,\om)$ contains left invertible, non-invertible elements.
\end{section}

\begin{section}{An application of direct finiteness}
We fix some notation. Throughout $\Gm$ will denote a discrete group $\Gm$, with identity element $e$. Given $p \in [1,\infty)$ we denote by $\lp{p}(\Gm)$ the Banach space of all $p$-summable functions $\Gm \to \Cplx$. $\lp{\infty}(\Gm)$ denotes the Banach space of all bounded functions $\Gm \to \Cplx$, equipped with the supremum norm. (The choice of complex scalars is not important: one could work throughout over $\Real$.)

Our main result is the following.
\begin{thm}\label{t:mainthm}
Let $f \in \lp{1}(\Gm)$ and let $\rho_f: \lp{\infty}(\Gm)\to\lp{\infty}(\Gm)$ denote the operator given by $h \mapsto h*f$. Then the following are equivalent:
\begin{enumerate}
\item $\rho_f$ is injective;
\item $\rho_f$ is bijective;
\item $f$ is invertible in $\lp{1}(\Gm)$, i.e.\ there exists $g \in \lp{1}(\Gm)$ such that $f*g=\delta_e=g*f$.
\end{enumerate}
\end{thm}
Theorem~\ref{t:mainthm} follows quickly from an old result, essentially due to Kaplansky, that the Banach algebra $\lp{1}(\Gm)$ is \dt{directly finite}. We recall that a unital ring $R$ is said to be \dt{directly finite} (or \dt{Dedekind finite} \cite{Lam}, or \dt{von Neumann finite}) if every left invertible element in $R$ is invertible: that is, if every pair $x, y \in R$ satisfying $xy=\id$ also satisfies $yx=\id$.

It is known \cite[p.\ 122]{Kap_FR} that  the group von Neumann algebra $VN(\Gm)$ -- and therefore any unital subring such as $\lp{1}(\Gamma)$ -- is directly finite. (A proof that $\lp{1}(\Gm)$ is directly finite, using purely $C^*$-algebraic methods, may be found in \cite{Mont}.) The key points are that: (i) in a $C^*$-algebra, every idempotent is similar to a hermitian idempotent; and (ii) there is a faithful tracial state on $VN(\Gm)$. See, e.g.~\cite[\S6.7]{KadRing_v2} for background and more details.

To apply Kaplansky's result we need the following elementary lemma. We shall denote the dual space of a Banach space $E$ by $E'$.

\begin{lemma}\label{l:dense-implies-onto}
Let $A$ be a unital Banach algebra and let $f \in A$. Let $R_f: A \to A$ be the operator on $A$ given by right multiplication by $f$, and let ${R_f}':A' \to A'$ be the adjoint operator. Then ${R_f}'$ is injective if and only if $R_f$ is surjective.
\end{lemma}
\begin{proof}
Sufficiency is immediate. To prove necessity, assume ${R_f}'$ is injective: then by the Hahn-Banach theorem $R_f$ has dense range. Since the group of units in $A$ is an open subset of $A$, the range of $R_f$ contains an invertible element, $u$ say. Thus there exists $h \in A$ with $hf=u$, and so for any $x \in A$ we have
\[ R_f (xu^{-1}h)= xu^{-1}hf = x \]
so that $R_f$ is surjective as claimed.
\end{proof}

%We pause to introduce some convenient notation. Let $\alpha:\Gm\to\Gm$ denote the map $t\mapsto t^{-1}$, and let $\alpha^*:\lp{1}(\Gm)\to\lp{1}(\Gm)$ be the corresponding involution on $\lp{1}(\Gm)$. A little thought shows that
%\[ \alpha^*(f*g) = \alpha^*(g)\alpha^*(f) \qquad\text{ for all $f,g\in\lp{1}(\Gm)$}. \]

\begin{proof}[Proof of Theorem~\ref{t:mainthm}]
The implications $(3)\implies(2)\implies(1)$ are obvious. We shall prove that $(1)\implies(3)$.

Thus, assume $\rho_f$ is injective. Let $R_f:\lp{1}(\Gm)\to \lp{1}(\Gm)$ denote the operator of right convolution with $f$, i.e.
\begin{equation}\label{eq:R-def}
R_f(h)(z) = (h*f)(z) = \sum_{y\in\Gm} h(zy^{-1})f(y) 
\quad\text{for $h\in\lp{1}(\Gm)$ and $z\in \Gm$.}
\end{equation}
The duality between $\lp{1}(\Gm)$ and $\lp{\infty}(\Gm)$ allows us to identify $\rho_f$ with ${R_f}'$ (this is easily verified by a direct comparison of \eqref{eq:rho-def} and \eqref{eq:R-def}). Hence by Lemma~\ref{l:dense-implies-onto}, $R_f$ is surjective. In particular there exists $g \in \lp{1}(\Gm)$ such that
\[ \delta_e = R_f(g) = g*f \quad; \]
and since $\lp{1}(\Gm)$ is directly finite, $g*f=\delta_e=f*g$, as required.
\end{proof}
\end{section}

\begin{section}{Two simple counterexamples}
Our first example shows that working on $\lp{\infty}$ is crucial.

\begin{eg}\label{eg:no-lp}
%[``Failure'' of Theorem~\ref{t:mainthm} for $\lp{p}$, $p<\infty$.]
Let $f=\delta_0-\delta_1 \in \lp{1}(\Zahl)$. It is clear that the kernel of $\rho_f:\lp{\infty}(\Zahl)\to\lp{\infty}(\Zahl)$ is spanned by the constant function $\Zahl\to \{1\}$, and so regarded as an operator on $c_0$ $\rho_f$ is injective. On the other hand, a direct computation shows that for any $a\in c_0(\Zahl)$ we have $\rho_f(a)\neq\delta_0$.

Thus, if $E=\lp{p}(\Zahl)$ ($1\leq p < \infty$) or $E=c_0(\Zahl)$ then the convolution operator $\rho_f:E\to E$ is injective but not surjective.
\end{eg}

For our second example we need some background. If $G$ is a locally compact group with right-invariant Haar measure $\mu$, and $f\in L^1(G)$, we may define\footnotemark\ $R_f:L^1(G)\to L^1(G)$ by
\footnotetext{This formula for convolution differs from the `usual' one because we are using a right-invariant Haar measure rather than the `usual' left invariant one. Our choice is made for ease of comparison with Equations \eqref{eq:rho-def} and \eqref{eq:R-def}, and ensures that $R_\bullet$ defines a right action of $L^1(G)$ on itself.}
\[ R_f (h)(z) = (h*f)(z) = \int_G h(zy^{-1})f(y)\,d\mu(y) \]
and take $\rho_f={R_f}': L^\infty(G) \to L^\infty(G)$, given by
\[ \rho_f(g)(x) = \int_G g(xy) f(y)\,d\mu(y) \]

\begin{eg}\label{eq:no-lc}
We work on $G=\Torus$ with its usual topology. Let $f$ be any element of $L^1(\Torus)$ such that
\[ \FC{f}(n)\neq 0 \text{ for all $n\in\Zahl$.} \]
Then $\rho_f:L^\infty(\Torus)\to L^\infty(\Torus)$ is injective but has non-closed range (and in particular is not surjective).
\end{eg}

This observation may be folklore; in any case we supply a short proof using some basic results from Fourier analysis on $\Torus$. First note that by standard duality arguments for Banach spaces (see, e.g.~\cite[Theorem 4.14]{RudFA}) it is enough to show that $R_f:L^1(\Torus)\to L^1(\Torus)$ has dense range but is not surjective. By taking Fourier transforms, this is equivalent to showing that the multiplication operator 
\[ M_\FC{f}: \widehat{L^1(\Torus)} \to \widehat{L^1(\Torus)} \]
has dense range but is not surjective.

To see that $M_\FC{f}$ has dense range, note first that the range clearly contains all finitely supported sequences on $\Zahl$ (because all the coefficients of $M_\FC{f}$ are nonzero); then note that the space of all such sequences is dense in $\widehat{L^1(\Torus)}$ (since the trigonometric polynomials are dense in $L^1(\Torus)$). Now suppose that $M_\FC{f}$ is surjective: then in particular there exists $h\in L^1(\Torus)$  with $M_\FC{f}(\FC{h})=\FC{f}$, i.e.~such that
\begin{equation}\label{eq:dagger}
\tag{$\dagger$} \FC{f}(n)\FC{h}(n) = \FC{f}(n)\quad\text{ for all $n\in\Zahl$.} \end{equation}
By the Riemann-Lebesgue lemma, since $h\in L^1(\Torus)$ the Fourier coefficients of $h$ tend to zero; in particular there is some $m$ for which $\abs{\FC{h}(m)} < 1$. Combined with \eqref{eq:dagger} this forces $\FC{f}(m)=0$, which contradicts our initial assumption on $f$; hence $M_\FC{f}$ is not surjective.

\end{section}

\begin{section}{Partial results for the weighted case}
Suppose we have a weight function $\om$ on $\Gm$, i.e.\ a function $\om: \Gm \to (0,\infty)$ such that $\om(xy)\leq\om(x)\om(y)$ for all $x,y \in \Gm$. Then given $f$ satisfying
\begin{equation}\label{eq:wt-l1}
 \sum_{x \in \Gm} \abs{f(x)} \om(x) < \infty
\end{equation}
one may define a weighted analogue of $\rho_f$, namely the operator
\begin{equation}\label{eq:wt-conv}
 \rho^\om_f : \lp{\infty}(\Gm) \to \lp{\infty}(\Gm) \end{equation}
given by 
\[ \rho^\om_f(g)_x =  \sum_{y \in \Gm} \frac{\om(xy)}{\om(x)} g(xy)f(y)  \quad\quad\left(g \in \lp{\infty}(\Gm);\; x \in \Gm\right) \]
(note that with this notation, when $\om$ is the constant weight 1 then $\rho^\om_f\equiv \rho_f$.)

It is natural to seek to extend Theorem~\ref{t:mainthm} to operators of the form $\rho^\om_f$. As before we may identify $\rho^\om_f$ with the adjoint of a multiplier $R_f$ on a certain unital Banach algebra, namely $\lp{1}(\Gm,\om)$.

More precisely, we have the following.
\begin{lemma}
Let $A$ denote the weighted convolution algebra $\lp{1}(\Gm,\om)$, and identify the dual of $A$ with the Banach space $\lp{\infty}(\Gm,\om^{-1})$. Let $M_\om:\lp{\infty}(\Gm)\to A'$ be the isometric isomorphism of Banach spaces defined by
\[ (M_\om h)_t=h(t)\om(t) \qquad(h\in\lp{\infty}(\Gm), t\in \Gm). \]

Let $f\in A$ and let $R_f:A\to A$ denote the operator of right multiplication by $f$.
Then $\rho^\om_f = M_\om^{-1} (R_f)' M_\om$.
\end{lemma}

The proof of this lemma consists of routine checking and is omitted.

To proceed as in the proof of Theorem~\ref{t:mainthm} we would need to know that $\lp{1}(\Gm,\om)$ is directly finite, and I do not know of any proof (or any counter\-example). However, for certain weights, or certain kinds of group, the arguments used to prove Theorem~\ref{t:mainthm} go through without problems.

To save needless repetition we introduce the foll\-owing terminology, following \cite{MCW}. If $\Gm$ is a discrete group and $\om$ is a weight function on $\Gm$, we refer to the pair $(\Gm,\om)$ as a \dt{weighted group}. We say that the weighted group $(\Gm,\om)$ is \dt{directly finite} if the convolution algebra $\lp{1}(\Gm,\om)$ is directly finite.

\begin{lemma}\label{l:wt-bounded-below}
Let $\om$ be a weight function on $\Gm$ such that
\begin{equation}\label{eq:wt-condn}
\inf_{x \in \Gamma} \om(x) > 0 
\end{equation}
(this happens, for example, if the weight $\om$ is symmetric). Then the weighted group $(\Gm,\om)$ is directly finite.
\end{lemma}
\begin{proof}
Condition \eqref{eq:wt-condn} ensures that $\lp{1}(\Gm,\om)$ \emph{is contained in} $\lp{1}(\Gm)$ as a unital subalgebra. As in the remarks following Theorem~\ref{t:mainthm} we deduce that $\lp{1}(\Gm,\om)$, being a unital subring of a directly finite ring, is itself directly finite.
\end{proof}

%As observed earlier, a weighted group will be directly finite if the weight $\om$ satisfies $\inf_{x\in\Gm}\om(x) > 0$. 
Using a result from \cite{MCW}, we can build on this lemma slightly, as follows.
\begin{propn}\label{p:wt-amen_implies_DF}
Let $(\Gm,\om)$ be a weighted group where $\Gm$ is amenable. Then $(\Gm,\om)$ is directly finite.
\end{propn}

\begin{proof}
By \cite[Lemma 1]{MCW} there exists a multiplicative function $\varphi:\Gm\to (0,\infty)$ such that $\varphi(x)\leq\om(x)$ for all $x\in\Gm$. Define a new weight $\nu$ on $\Gm$ by setting $\nu(x)=\varphi(x)^{-1}\om(x)$ for all $x$; then clearly $\inf_{x\in\Gm} \nu(x) \geq 1$, and so by Lemma~\ref{l:wt-bounded-below} the weighted group $(\Gm,\nu)$ is directly finite. Now observe that there is a continuous algebra isomorphism $\theta$ from $\lp{1}(\Gm,\om)$ onto $\lp{1}(\Gm,\nu)$, given by
\[ (\theta a)(x) = \varphi(x)a(x) \qquad(a\in\lp{1}(\Gm,\om)\/;\/ x\in\Gm) \]
and therefore $\lp{1}(\Gm,\om)$, being isomorphic to a directly finite algebra, is itself directly finite.
\end{proof}

\begin{thm}
Suppose that the weighted group $(\Gm,\om)$ is directly finite.
Let $f$ and $\rho^\om_f$ be as given in \eqref{eq:wt-l1}, \eqref{eq:wt-conv}. Then the following are equivalent:
\begin{itemize}
\item[{\rm (i)}]   $\rho^\om_f$ is injective;
\item[{\rm (ii)}]  $\rho^\om_f$ is bijective;
\item[{\rm (iii)}] $f$ is invertible in $\lp{1}(\Gm.\om)$, i.e.\ there exists $g:\Gm\to\Cplx$ with
\[  \sum_{x \in \Gm} \abs{g(x)} \om(x) < \infty \]
such that $f*g=\delta_e=g*f$.
\end{itemize}
\end{thm}
\begin{proof}[Outline of proof]
We only need to prove that (i)$\implies$ (iii). As before, we use Lemma~\ref{l:dense-implies-onto} to show that if $\rho_f$ is injective then there exists $g \in \lp{1}(\Gm,\om)$ such that $f*g=\delta_e$. Since $(\Gm,\om)$ is directly finite we deduce that $g*f=\delta_e$ as required.
\end{proof}

\begin{rem}
We do not know if every weighted group is directly finite. In the proof of Proposition~\ref{p:wt-amen_implies_DF} we relied crucially on the fact that any weight on an amenable group dominates a non-zero character on that group. In \cite[Example 2]{MCW} an example is given of a weighted group for which the only character dominated by the weight is zero, and so any attempt to prove that all weighted groups are directly finite must use deeper arguments than are employed here. 
\end{rem}
\end{section}

\subsection*{Acknowledgements} I would like to thank C. Deninger, G. Elek and N.~J. Young for helpful exchanges, and in particular F. Ghahramani for reminding me of the main result in~\cite{MCW}.

% bibliography cut and pasted from bibtex output

\end{document}